\documentclass[11pt,reqno]{amsart}

\usepackage{amsmath, amsthm, amssymb,amsaddr}
\usepackage{dcolumn}
\usepackage{bm}
\usepackage{amsmath} 
\usepackage{mathrsfs}
\usepackage{mathtools}
\usepackage{algpseudocode}
\usepackage{algorithm}
\usepackage[pdftex]{graphicx}
\usepackage{subfigure}
\usepackage{epstopdf}
\usepackage{tikz}
\usepackage{xcolor}
\usepackage[export]{adjustbox}
\usepackage{caption}
\usepackage{hyperref}
\usepackage{stmaryrd}
\usepackage[utf8]{inputenc} 
\usepackage[T1]{fontenc}
\usepackage{enumitem}

\usepackage{titlesec}

\titleformat{\section}{\bfseries\Large}{\thesection}{1em}{}
\titleformat{\subsection}{\bfseries\large}{\thesubsection}{1em}{}
\titleformat{\subsubsection}{\bfseries\normalsize}{\thesubsubsection}{1em}{}

\theoremstyle{remark}

\newtheoremstyle{italic_lemma} 
{} 
{} 
{\itshape} 
{} 
{\itshape} 
{.} 
{ } 
{} 

\theoremstyle{italic_lemma}

\newcommand{\be}[1]{\begin{equation}\label{#1}}
\newcommand{\ee}{\end{equation}}

\newcommand{\no}[1]{#1}

\newcommand{\bx}{\mathbf{x}}

\newcommand{\bk}{\boldsymbol{k}}

\newcommand{\iu}{{i\mkern1mu}}

\renewcommand{\no}[1]{} 

\newcounter{subsubsubsection}[subsubsection]

\title[A MATLAB Package for Quantitative Ultrasound Tomography via Ray-Born Inversion]{Introduction and Numerical Validation of an Open-Source MATLAB Package for Quantitative Ultrasound Tomography via Ray-Born Inversion}
\author{Ashkan Javaherian}
\email{ajavaherian62@gmail.com; ashkan.javaherian@ut.ac.ir}
\address{
Department of Bio-Electric, School of Electrical and Computer Engineering,
University College of Engineering, University of Tehran, Tehran, Iran.
\footnote{The \textsc{MATLAB} codes supporting the numerical results presented in this study are fully reproducible by running the example script \texttt{example\_validate\_ray\_fish\_eye\_phantom.m} available in the open-source GitHub repository:
\url{https://github.com/Ash1362/ray-based-quantitative-ultrasound-tomography}.}
}
\date{November 2025}

\begin{document}
\maketitle

\section*{Abstract}
We present a \textsc{MATLAB} package for reconstructing sound-speed images from transmission ultrasound data. The package is based on two-point ray tracing and implements two complementary inversion strategies for image reconstruction. The first is a time-of-flight (ToF) method that produces low-resolution, low-contrast images with minimal artefacts. The second is a ray–Born inversion method, which integrates high-frequency ray theory with the Born approximation to generate high-resolution sound-speed reconstructions. Early iterations of the ToF reconstruction are used to provide an initial estimate for the more advanced ray–Born approach. The core of this software package consists of four ray-tracing algorithms, whose accuracy is assessed in this study with respect to known analytical trajectories and accumulated acoustic path lengths. Furthermore, both image-reconstruction strategies have been validated numerically with simulated synthetic datasets and experimentally with open-source in-vitro and in-vivo datasets in related parallel studies.
\noindent

\section{Introduction}
\noindent

Ultrasound tomography aims to reconstruct maps of acoustic properties inside an object using pressure data measured in time by ultrasonic transducers placed around the object \cite{Li-tof-imag}. The imaging process involves a set of transducers placed around the region of interest. Acoustic waves are generated by exciting one transducer at a time; these waves propagate through the object and are recorded by the remaining transducers. This excitation and measurement cycle is repeated successively for all transducers in the array \cite{Gemmeke}.

The acoustic maps reconstructed from ultrasound data can be either \textit{qualitative} or \textit{quantitative}. Qualitative imaging typically identifies inter-medium boundaries based on echoes reflected from interfaces where acoustic impedance changes \cite{Synnevag}. In contrast, quantitative imaging seeks to estimate intrinsic material properties, such as the spatial distribution of acoustic absorptivity or the speed of sound, quantitatively \cite{Wiskin-nonlinear}. These properties are primarily inferred from the amplitude and phase of the transmitted pressure waves measured by the transducers.

This manuscript introduces an open-source MATLAB package for ultrasound tomography \cite{Javaherian0,Javaherian00} using two-point ray tracing \cite{Cerveny}. While the focus is on the toolbox's application to quantitative reconstruction of the speed of sound, it can be extended to support other ultrasound imaging modalities. The image reconstruction approaches included in this toolbox have been validated using simulated transmission datasets \cite{Javaherian-refraction,Javaherian-scattering,Javaherian-scattering2}, along with in-vitro and in-vivo datasets released by the University of Rochester Medical Center \cite{Javaherian-experimental}.

Section 2 outlines the image reconstruction algorithms for generating quantitative images from transmission ultrasound data. Section 3 introduces the presented \textsc{MATLAB} software package and describes its various processing steps. Section 4 evaluates the accuracy of the four ray-tracing algorithms against analytic ray trajectories and accumulated acoustic path lengths. Section 5 summarizes the significance of the developed software package, while the underlying theory is briefly outlined in the Appendix.

\section{Image Reconstruction Approaches}

\noindent
Approaches for sound speed image reconstruction from ultrasound data are generally categorized based on the type of information extracted from the measured time series.

\subsection{Time-of-Flight (ToF)–Based Methods}

These traditional approaches rely solely on the first arrivals of the measured time series, which correspond to the shortest time it takes for acoustic waves to travel between each emitter–receiver pair—also known as the time of flight (ToF). These methods are robust and computationally inexpensive but generally produce low-resolution images of the sound speed distribution \cite{Li-tof-imag,Li-tof-pick,Qu,Ceccato}.

Time-of-flight (ToF)-based approaches typically aim to reconstruct a sound speed map by iteratively updating the model until the misfit between the measured and modeled ToFs is minimized in the least-squares (\(\ell_2\)-norm) sense. The measured ToFs are extracted once, prior to inversion, by identifying the first arrival times in the time-domain pressure traces. The first-arrival picking is performed separately for each emitter–receiver pair \cite{Li-tof-pick}.

Each modeled time-of-flight (ToF) corresponds to the accumulated travel time along the ray connecting a specific emitter–receiver pair. These rays can be assumed to be straight \cite{Dapp} or bent \cite{Li-tof-imag,Huang}, depending on the modeling assumptions. In the latter case, ray paths bend due to refraction caused by heterogeneities in the refractive index (or wavenumber) distribution, typically under a high-frequency approximation \cite{Javaherian-refraction}.

\subsection{Full-Waveform Inversion (FWI)} These approaches utilize the entire recorded time traces for image reconstruction. FWI methods can produce high-resolution images but significantly increase computational expense due to the need to model full wave physics and the complexity of the associated optimization procedures \cite{Ali-2024,Robins}.

\subsection{Singly-Scattering–Based Methods}

The third class of approaches for reconstructing the sound speed map from ultrasound data incorporates the primary scattered waves into image reconstruction. Instead of solving the full nonlinear wave equation, these methods approximate the scattered field as a first-order perturbation around a known reference medium \cite{Simonetti}. This linearization implies that the scattered wavefield depends linearly on the contrast (difference) between the true medium and the reference. The linearization is accurate only if the perturbations in sound speed are sufficiently small and higher-order scattering (multiple scattering) is negligible \cite{Huthwaite}. At higher frequencies, the wavelength shortens, causing waves to interact more strongly and repeatedly with inhomogeneities in the medium, thereby producing multiple scattering. This effect violates the single scattering assumption underlying the linearization at higher frequencies.


\section{Software Description}  \label{sec:software-description}

This \textsc{MATLAB} package is primarily designed for imaging weakly heterogeneous media, such as soft tissues. This toolbox implements two classes of methods for reconstructing sound speed images from ultrasound datasets. 

\begin{enumerate}
    \item \textit{ToF-Based Inversion:}  
    The initial aim of developing this software package was to extend traditional time-of-flight (ToF)–based reconstruction methods—originally developed for 2D ultrasound data acquired using a ring-shaped array \cite{Ali-2024} or a rotating pair of opposite linear arrays \cite{Robins}—to 3D imaging scenarios \cite{Javaherian-refraction}. These extensions enable the fast reconstruction of low-resolution sound speed images using more complex acquisition geometries, such as hemispherical arrays \cite{Gemmeke,Javaherian-refraction} or rotating pairs of opposite planar arrays \cite{Wiskin-nonlinear}. 

    \item \textit{Ray-Born Inversion:}  
    The toolbox was later enhanced to support high-resolution image reconstruction by combining ray theory with the Born approximation \cite{Javaherian-scattering,Javaherian-scattering2,Javaherian-experimental}. The frequency range typically used in ultrasound tomography for biomedical applications spans approximately 0.3 to 1.2 MHz, covering both relatively low and moderately high frequencies \cite{Javaherian-experimental,Ali-2024}. Motivated by this, we developed image reconstruction approaches that combine the Born approximation with ray theory based on high-frequency asymptotics \cite{Javaherian-scattering,Javaherian-scattering2,Thierry}. Unlike the standard Born approximation, which assumes Green’s functions in a homogeneous water background \cite{Simonetti}, our method approximates them along bent rays. This accounts for phase aberrations due to acoustic refraction and dispersion, as well as amplitude decay caused by geometrical spreading and acoustic absorption \cite{Javaherian-scattering,Javaherian-scattering2}. Consequently, the medium’s heterogeneity affects not only the scattering potential but also the Green’s functions, leading to linearized equations that more accurately capture the physics of wave propagation. Since the Born approximation is valid only for small perturbations in sound speed, an initial estimate is first obtained using a ToF-based method \cite{Javaherian-scattering2,Javaherian-experimental}. Our numerical experiments \cite{Javaherian-scattering,Javaherian-scattering2}, as well as implementations on in vitro and in vivo datasets \cite{Javaherian-experimental}, demonstrate that the early iterations of the ToF-based inversion provide a sufficiently accurate initial estimate for the subsequent ray–Born reconstruction.
\end{enumerate}


In contrast to full-wave inversion approaches—which often function as black-box procedures for solving the wave equation and the associated inverse problem \cite{Ali-2024}—our developed ray-based image reconstruction methods consist of a sequence of well-defined algorithmic steps. While this modular structure enhances flexibility, it may also appear complex to new users. This section outlines the key components implemented in the software to facilitate the reconstruction process. For detailed descriptions of each step, users are referred to the relevant publications by the author.

\subsection{Minimization and Linearizations} \label{sec:minimization}

Below, the minimization processes for the ToF-based and Ray-Born image reconstruction approaches are described.

\subsubsection{Linearizations in ToF-based Inversion}
In time-of-flight (ToF) based approaches, the speed of sound is iteratively updated until a misfit function, defined as the \( \ell_2 \)-norm of the discrepancy between the measured and modeled time-of-flights (ToFs), is minimized \cite{Javaherian-refraction}. Since the forward operator that maps the sound speed to the modeled ToFs is generally nonlinear, the minimization problem must be iteratively linearized. At each iteration, solving the associated linearized subproblem yields an update direction, which can be computed using different numerical strategies.

These numerical algorithms are inherently iterative, requiring inner recursions to solve each linearized subproblem. Each linearization involves computing the Jacobian matrix. The Jacobian matrix—also referred to as the \emph{system matrix}—maps perturbations in the sound speed distribution at discrete grid points to the corresponding perturbations in the modeled times-of-flight (ToFs) for all emitter–receiver pairs. It is a sparse matrix of size \( N_e N_r \times N_n \), where \( N_e \), \( N_r \), and \( N_n \) denote the number of emitters, receivers, and grid points, respectively. Because of its sparsity, it can be computed and stored explicitly. 
Each row of the system matrix represents the contribution of the grid points to the acoustic path length of the ray connecting a specific emitter–receiver pair \cite{Javaherian-refraction}. Because ray tracing is performed off-grid, the contributions of nearby grid points are weighted using interpolation coefficients corresponding to the sampled points along each ray, while the contributions of grid points far from the ray are set to zero \cite{Javaherian-refraction,Anderson}. 

Our software package implements two numerical approaches for solving these linearized subproblems:

\begin{enumerate}
    \item \textit{Conjugate Gradient (CG) Method:} CG is a widely used iterative technique for solving large-scale linear minimization problems. In our application, it often reconstructs sharper contrast and reveals more features in the updated sound speed maps. However, it may introduce artifacts due to neglecting spatial heterogeneity in ray density. This may result in enhanced contrast in regions with dense ray coverage, but reduced contrast in areas with sparse coverage, including acoustic shadow regions.
    
    \item \textit{Simultaneous Algebraic Reconstruction Technique (SART):} Specifically developed for ray-based tomography, SART compensates for nonuniform ray sampling by applying appropriate weighting to the system matrix in both the data and solution spaces \cite{Anderson}. This balances contributions from rays across the imaging domain, resulting in smoother update directions for the sound speed with fewer artifacts. Analogous to the conjugate gradient (CG) approach, the linearized subproblem is subsequently minimized in an iterative manner. \textit{(Default)}
\end{enumerate}


\subsubsection{Linearizations in Ray-Born Inversion}
Image reconstruction is performed in the frequency domain. Starting from low frequencies, the misfit function—defined as the $\ell_2$-norm of the discrepancy between the measured and modeled data—is linearized, and the resulting linear equation is solved. The updated sound-speed map at each frequency level is then used to solve the linearized subproblem at the next higher frequency. For each linearized subproblem (frequency level), the ray trajectories are updated according to the current estimate of the wavenumber field. By progressing from low to high frequencies, the magnitude of the sound-speed update (in the $\ell_2$-norm sense) gradually decreases until it falls below a predefined threshold~\cite{Javaherian-scattering,Javaherian-scattering2,Javaherian-experimental}.

Solving each linearized subproblem is equivalent to computing the action of the inverse Hessian matrix on the negative gradient of the objective function. Our toolbox implements two numerical approaches for this purpose:

\begin{enumerate}

\item \textit{Hessian-Based (Iterative):}  
At each frequency level, the linearized equation is solved iteratively using a linear Conjugate Gradient (CG) algorithm with multiple inner iterations. Unlike the ToF-based method, the Jacobian matrix in this approach is not sparse and must therefore be computed implicitly. The computational cost of each inner iteration is comparable to that of the ray tracing performed at the beginning of each linearization. Consequently, the overall cost increases by approximately an order of magnitude~\cite{Javaherian-scattering}.

\item \textit{Hessian-free (single-step):} Rather than performing an implicit inversion of the Hessian matrix through iterative methods, this approach applies weighting to the linear equation in both data and solution spaces, rendering the Hessian matrix approximately diagonal and invertible in a single step \cite{Javaherian-scattering2,Javaherian-experimental}. Consequently, the primary computational expense at each frequency level corresponds to ray tracing, reducing the total cost by roughly an order of magnitude compared to the Hessian-based method \cite{Javaherian-scattering2}. (\textit{Default})

\end{enumerate}

\subsection{Ray Tracing} \label{sec:ray-tracing}

In both the ToF-based and ray-Born image reconstruction approaches, the forward operator is linearized and modeled using ray tracing. The details of the ray tracing procedure are described below.

\subsubsection{Ray Tracing in ToF-based Inversion}
At each linearized subproblem, ray tracing is performed to compute the accumulated travel time between each emitter-receiver pair. Additionally, the contribution of each grid point to this travel time is weighted by the ray-to-grid interpolation coefficients \cite{Javaherian-refraction}.

\subsubsection{Ray Tracing in Ray-Born Inversion}  
At each frequency level, the corresponding linearized subproblem is formulated and solved using forward and backward Green's functions, which are approximated along rays initialized from all emitter and receiver positions, respectively \cite{Javaherian-scattering, Javaherian-scattering2}. These rays are used to approximate the Green's function parameters throughout the domain, including:
\begin{enumerate}
    \item the accumulated phase, accounting for distortions due to dispersion and caustic effects; and
    \item the amplitude, incorporating geometrical spreading and accumulated acoustic absorption.
\end{enumerate}

\subsubsection{Ray Tracing Approaches}  
In prototype ToF-based implementations, ray trajectories are computed on a discrete grid using grid-based ray tracing techniques \cite{Li-tof-imag,Klimes}. While such approaches are straightforward to implement, they suffer from two important limitations. First, the optimal ray paths—those that minimize acoustic path length—are not necessarily confined to trajectories that follow connected grid points. This mismatch can introduce errors in both forward modeling and inversion, often requiring the use of finer grids with reduced spacing to improve accuracy, which in turn increases computational cost. The second limitation arises in sparse transducer arrays, where rays initialized from each emitter are traced across all grid points, even though the system matrix is ultimately constructed only from rays that connect emitters to receivers. This unnecessary tracing adds computational burden without contributing directly to the image reconstruction. On the other hand, grid-based ray tracing offers improved accuracy near high-contrast regions and sharp interfaces, because ray bending follows fixed rules on the grid and is not affected by divergence instabilities present in off-grid formulations.


For weakly heterogeneous media, such as soft tissues, the trajectory of a ray connecting an emitter–receiver pair can be determined using two-point ray tracing \cite{Cerveny,Javaherian-refraction}. This corresponds to a boundary value problem, in which the ray path is iteratively updated until a termination criterion is satisfied. Two common classes of two-point ray tracing methods are the \textit{bending method} \cite{Moser} and the \textit{shooting method} \cite{Rawlinson}. In the bending method, the endpoints of the ray are fixed to the emitter and receiver positions, which are typically treated as points. The ray trajectory is iteratively adjusted to minimize the acoustic path length between these points \cite{Moser}. A major challenge with this approach is defining a reliable termination condition and ensuring convergence of the associated iterative algorithm. In contrast, the shooting method fixes the ray's starting position at the emitter and iteratively traces the ray by updating its initial propagation direction until it intercepts the receiver \cite{Rawlinson}. 

In our software package, the shooting method is adopted owing to its flexibility, straightforward implementation, and robust convergence in weakly heterogeneous media \cite{Javaherian-refraction}. Using shooting methods, each update of a ray's trajectory amounts to solving an initial value problem, given the ray's initial position and a unit direction vector \cite{Javaherian-refraction,Rawlinson}. In this context, ray tracing refers to computing the path of a ray that minimizes the acoustic path length, in accordance with Fermat's principle. This principle underlies the high-frequency approximation of wave propagation, leading to the Eikonal equation as the governing model \cite{Cerveny,Anderson2}. The resulting ray equations compute the evolution of the position \(\mathbf{x}\) and unit direction \(\mathbf{d}\) along the ray, depending explicitly on the spatial distribution of the refractive index (or equivalently, the wavenumber) and its gradient.

Our software package supports four numerical algorithms for ray tracing:
\begin{enumerate}
    \item The ``Dual-Update'' method (see \cite{Javaherian-refraction}, Section 3.3.1),
    \item The ``Mixed-Step'' method (see \cite{Javaherian-refraction}, Section 3.3.2, or \cite{Anderson2}, Section II-A),
    \item The Method of ``Characteristics'' (see \cite{Anderson2}, Section II-B), and
    \item The ``2nd order Runge--Kutta (Heun's)'' method (see \cite{Javaherian-scattering}, Section 5.2), (\textit{Default}).
\end{enumerate}

Using the ToF-based image reconstruction approach, Heun's method is employed as the default, although all the listed ray tracing methods are supported. For the Ray-Born image reconstruction approach, the software exclusively utilizes Heun's method.

\subsection{Grid-to-Ray Interpolation}

While the sound speed and its gradient are defined and updated on discrete grid points, the sampled points along each ray—computed using the ray tracing algorithms discussed above—are generally off-grid and can lie anywhere within the domain. 

Furthermore, in the \textit{ToF-based} image reconstruction approach, each row of the Jacobian matrix corresponds to an emitter–receiver pair and is weighted by interpolation coefficients associated with grid points near the corresponding ray path.

\subsubsection{Interpolation Approaches}

Our software package supports two interpolation schemes:

\begin{enumerate}
    \item \textit{Bilinear interpolation:} This method is computationally efficient. However, the gradient of the interpolated field is obtained by directly interpolating values defined on discrete grid points, which introduces considerable errors. For each off-grid point, both the field and its gradient are interpolated using only the vertices of the voxel enclosing the point—sacrificing accuracy for computational efficiency (see \cite{Javaherian-refraction}, Section 3.4.). 
    
    \item \textit{Cubic B-spline interpolation:} By expressing the field as a piecewise third-order polynomial defined by control points located at the grid nodes, cubic B-spline interpolation produces a \( \mathcal{C}^2 \)-continuous interpolated field. Both first- and second-order derivatives are computed analytically. For each off-grid point, interpolation involves a broader stencil that includes four neighboring grid points. While computationally more expensive, this approach provides significantly improved accuracy, particularly for gradients (see \cite{Javaherian-scattering}, Section 5.3.).
\end{enumerate}

\subsubsection{Interpolation in ToF-based Inversion}

Since ToF-based methods inherently produce low-resolution images of the sound speed, \textit{Bilinear} interpolation is used by default in our toolbox. This choice balances accuracy and computational efficiency \cite{Javaherian-refraction}.

\subsubsection{Interpolation in Ray-Born Inversion}

For Ray–Born image reconstruction, a paraxial ray-tracing system must be solved, in which the wavenumber and its first and second gradients are evaluated at sampled points along each ray. To achieve high accuracy, the toolbox employs a \textit{cubic B-spline} interpolation approach, which computes analytic gradients—essential for precise ray tracing and amplitude approximations \cite{Javaherian-scattering, Javaherian-scattering2,Javaherian-experimental}.

\subsection{Ray Linking}
As explained in Section \ref{sec:ray-tracing}, determining the trajectory of a ray connecting each emitter-receiver pair is a boundary value problem, which is here solved using shooting methods, also known as ray linking \cite{Javaherian-refraction}. For each emitter-receiver pair, a misfit function is defined in terms of the discrepancy of the interception point of the ray by the boundary after traveling across the object and the receiver's center. Ray linking is a root-finding problem, which is parameterized in the angular domain \cite{Javaherian-refraction}. For solving this problem, the initial unit direction of the ray in the angular domain is iteratively adjusted until this misfit function is zeroed in the angular domain \cite{Javaherian-refraction}. 

\noindent
For 2D configurations, the common ray-linking approaches include:

\begin{enumerate}
    \item \textit{Regula Falsi (False Position) Method:} 
    Two initial unit directions are chosen such that the optimal initial unit direction lies between them. Starting from these initial guesses, the interval is recursively reduced by replacing one of the directions with a new root estimate computed via linear interpolation. This process continues until the size of the interval falls below a predetermined tolerance level. This approach is robust to poor initial guesses and does not diverge for sufficiently smooth media. However, its convergence is relatively slow.

    \item \textit{Secant Method:} 
    The initial unit direction is iteratively updated using secant lines computed from the two most recent estimates. Each Secant direction approximates the inverse of a first-order gradient, computed from the two last updates of the initial unit direction using finite difference approximation, applied to the negative misfit function. This method generally converges faster than Regula Falsi but can be more sensitive to initial guesses. (\textit{Default})
\end{enumerate}

\noindent    
For 3D configurations, we proposed: 

\begin{enumerate}[label={}]
    \item 
\textit{Broyden-like Method:} This is a derivative-free method for solving nonlinear root-finding problems. In the \textit{Newton-Raphson method}, each update of the unknown parameter requires an explicit computation of the Jacobian matrix using finite differences. Therefore, for each update of the initial unit direction (in 3D), both the azimuthal and polar initial angles must be perturbed, requiring the tracing of at least three rays per update \cite{Javaherian-refraction}. 

In contrast, our proposed Broyden-like approach avoids explicit computation of the Jacobian matrix. Instead, the Jacobian is estimated and recursively updated using the two most recent updates of the initial unit directions, which reduces the computational cost to tracing only one ray per update \cite{Javaherian-refraction}.
\end{enumerate}

\subsubsection{Ray Linking in ToF-based Inversion}
For our application—quantitative reconstruction of the sound speed—we begin with a homogeneous initial guess for the sound speed and iteratively update it by solving a sequence of linearized subproblems, as described in Section~\ref{sec:minimization}. For each linearized subproblem, the ray linking problem must be solved for all emitter-receiver pairs separately. In the first subproblem, since the sound speed is assumed homogeneous, rays are initialized as straight lines directed toward the centers of the receivers, without requiring recursive ray linking. For each subsequent subproblem, the initial unit directions obtained through ray linking in the previous subproblem are used as initial guesses for the ray linking process \cite{Javaherian-refraction}. 

\subsubsection{Ray Linking in Ray-Born Inversion}
The ray linking procedure used in Ray-Born inversion follows the same strategy as in the Time-of-Flight (ToF)–based inversion. At each frequency level (i.e., for each linearized subproblem), ray linking is initialized using the unit direction vectors obtained from the previous frequency level. Rather than starting from a homogeneous background, the initial ray directions for the first frequency level are taken from the last linearization of the ToF-based method~\cite{Javaherian-scattering,Javaherian-scattering2,Javaherian-experimental}.


\subsection{Green's Function Approximation in Ray-Born Inversion}

In ray-Born image reconstruction methods, the linearized equation arising at each frequency level is expressed in terms of the forward and backward Green's functions, which are initialized at the emitter and receiver positions, respectively. The parameters of the Green’s functions, specifically the phase and amplitude, are approximated along ray paths.

Based on the \textit{Eikonal} equation, the accumulated phase is computed by integrating along the rays, with corrections for phase distortions caused by dispersion and caustics. The amplitude loss includes contributions from both geometric spreading and acoustic absorption along the propagation path. While acoustic absorption is accumulated along the linked (reference) rays, geometric spreading depends on the ray's Jacobian, as governed by the \textit{Transport} equation. 

The computation of the ray's Jacobian relies on information from rays in the vicinity of the reference ray. In our software package, two approaches are implemented for computing the ray's Jacobian:

\begin{enumerate}
\item \textit{Auxiliary rays:} Given the initial position of a ray at the emitter (or receiver) and its initial unit direction—previously computed in the angular domain through ray linking—auxiliary rays (two in the 2D case) are traced independently of the reference ray by perturbing the initial angle. The ray's Jacobian along the reference ray is then approximated from the trajectories of these auxiliary rays \cite{Cerveny,Javaherian-scattering}.

\item \textit{Paraxial ray:} Instead of tracing auxiliary rays independently of the reference ray, a paraxial ray is traced simultaneously with the reference ray by incorporating first-order perturbations into the ray equations, induced by small perturbations in the unit direction. The ray's Jacobian can then be directly computed from the solutions of the paraxial system, which are the perturbed positions \(\delta \mathbf{x}\) and unit directions \(\delta \mathbf{d}\). By leveraging information along the reference ray, paraxial ray tracing offers improved accuracy and computational efficiency compared to the auxiliary-ray approach~\cite{Cerveny,Javaherian-scattering2}. (\textit{Default})

\end{enumerate}

\begin{figure} 
\centering
\includegraphics[width=0.45\textwidth]{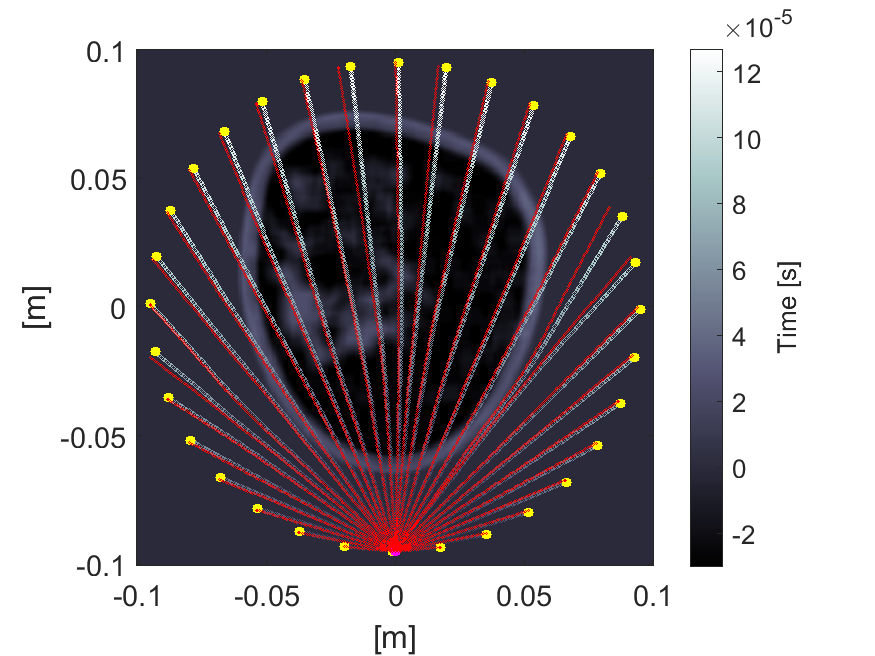}
 \caption{Trajectories of rays linking a single emitting transducer to all receiving transducers. The color of each linked ray is scaled by the accumulated time delay referenced to the emitter. Given the positions \(\mathbf{x}\) of the linked rays, the red-colored rays represent \(\mathbf{x} + 0.01\, \delta \mathbf{x}\), where \(\delta \mathbf{x}\) is computed by the paraxial ray tracing.}
\label{fig:1}
\end{figure}

Figure~\ref{fig:1} illustrates the rays connecting a transducer acting as an emitter to all other transducers acting as receivers. 
The rays are refracted due to spatial variations in the sound speed. In this figure, the color of each linked (reference) ray is scaled according to its accumulated time delay relative to the emitter. For each emitter–receiver pair, the initial direction of the ray is iteratively updated, and the ray is traced through a root-finding algorithm until its interception point coincides with the receiver.  For the ray that satisfies the termination criterion of this root-finding problem (i.e., the linked ray), a paraxial system of ray equations is solved. Given the ray position \(\mathbf{x}\) and unit direction \(\mathbf{d}\), the paraxial ray describes the first-order perturbation to the position and unit direction of the linked ray, denoted by \(\delta \mathbf{x}\) and \(\delta \mathbf{d}\), respectively. In Figure~\ref{fig:1}, the red-colored rays represent \(\mathbf{x} + 0.01\, \delta \mathbf{x}\). These rays are shown only for visualization purposes. In practice, the ray Jacobian is computed in terms of \(\delta \mathbf{x}\) and \(\delta \mathbf{d}\).

\section{Numerical Validation}

The numerical results and reconstructed images using both simulated and experimental data are presented in the author's other studies \cite{Javaherian-refraction,Javaherian-scattering,Javaherian-scattering2,Javaherian-experimental}. This section is dedicated to numerically validating the accuracy of the ray tracing algorithms, which constitute the core computational step in our image reconstruction approaches. To this end, we employ the \textit{Maxwell’s fish‐eye lens} phantom.

\subsection*{Maxwell's Fish-Eye Lens}

The accuracy of the ray tracing algorithms is assessed using the \textit{Maxwell's fish-eye lens} phantom, for which the ray trajectories are known analytically \cite{Anderson2}. Let \(\mathrm{dim} \) denote the number of spatial dimensions (i.e., \(\mathrm{dim} = 2\) or \(3\)). The position vector is defined as \(\mathbf{x} = [x_1, \dots, x_\mathrm{dim}]\), where \(x_1, \dots, x_{\mathrm{dim}}\) are the Cartesian coordinates. The Euclidean distance from the origin \(\mathbf{x}_o\) is denoted by
\[
r = \left\| \mathbf{x} \right\| = \sqrt{x_1^2 + \dots + x_{\mathrm{dim}}^2}.
\]

\noindent
The refractive index distribution of Maxwell's fish-eye lens is given by:
\begin{align}
    n(\mathbf{x}) = \frac{n(\mathbf{x}_o)}{1 + \left( \frac{r}{a} \right)^2},
\end{align}
where \(n(\mathbf{x}_o)\) is the refractive index at the origin, and is set to 1. The scalar parameter \(a\) defines the scale of the lens and is also set to 1 throughout this study. This rafractive index distribution is shown in figure~\ref{fig:2a}.

This phantom exhibits a key analytical property: ray trajectories follow perfect circular arcs that are tangent to their initial directions. A ray launched from a point \( \mathbf{x}_p \) propagates along a circular path, passing through a conjugate (mirror) point  \(\mathbf{x}_{\tilde{p}}\) before returning to its origin upon completing a full rotation. The center of this circular path, which varies with the ray’s initial direction, is denoted by \(\mathbf{x}_c\).

The computational grid is discretized with an angular resolution of 1 degree along the periphery of a circle with radius \(a\), resulting in a spatial grid spacing of
\[
\Delta x = \frac{2 \pi a}{360}.
\]

For the numerical results presented below, the refractive index distribution is defined on the grid, consistent with the procedures used in both the ToF-based and ray–Born image reconstruction algorithms implemented in the toolbox. Consequently, the errors (deviations) reported in this section include contributions from both the ray-tracing algorithm and the interpolation steps. In principle, these errors could be reduced to contributions solely from the ray-tracing algorithm by tracing rays in a Maxwell fish-eye lens defined in a continuous domain, for which the medium gradients can be computed analytically. However, such an evaluation appears unnecessary, since the discretized refractive index (or wavenumber) field is ultimately the one employed in practical image reconstruction frameworks. In the results presented below, both the grid-to-ray and ray-to-grid interpolations are carried out using cubic B-splines.

\subsection{Mean Radius Deviation}

One criterion for assessing the accuracy of traced rays is the \emph{radius deviation} of the circular paths they traverse. The ray’s initial position is defined as:
\begin{itemize}
    \item[1)] $\mathbf{x}_p = [0, a]$, \quad for \(\mathrm{dim}=2\),
    \item[2)] $\mathbf{x}_p = [0, 0, a]$, \quad for \(\mathrm{dim}=3\),
\end{itemize}
and the intended center of the ray’s circular path is:
\begin{itemize}
    \item[1)] $\mathbf{x}_c = [a, 0]$, \quad for \(\mathrm{dim}=2\),
    \item[2)] $\mathbf{x}_c = [a, a, 0]$, \quad for \(\mathrm{dim}=3\).
\end{itemize}

\noindent
Correspondingly, the initial direction of the rays is set as follows:
\begin{itemize}
    \item[1)] In 2D: a single ray with an initial unit direction vector oriented clockwise and normal to the axis \(\mathbf{x}_c - \mathbf{x}_p\), i.e., \(\theta = \pi/2\), where \(\theta\) is the angle between the initial direction and the axis.
    \item[2)] In 3D: 21 rays with initial unit directions normal to $\mathbf{x}_c - \mathbf{x}_p$, uniformly distributed over a full $2\pi$ rotation around the axis defined by $\mathbf{x}_c - \mathbf{x}_p$. The polar angle is fixed at $\theta = \pi/2$, while the azimuthal angles $\varphi$ are uniformly sampled over \([0, 2\pi)\).
\end{itemize}

\begin{figure} 
\centering
\subfigure[]{\includegraphics[width=0.40\textwidth]{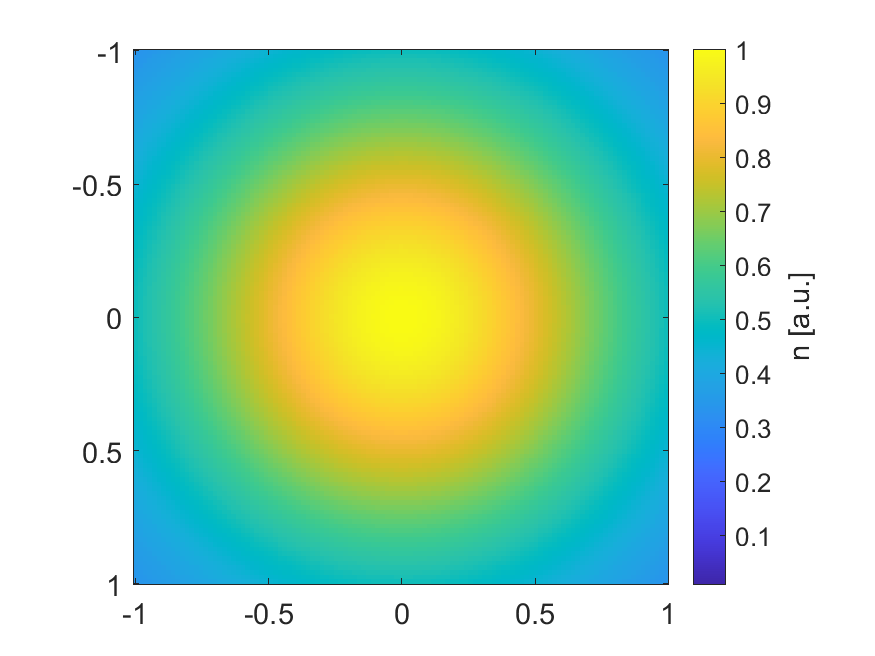}
\label{fig:2a}  }\\
\subfigure[]{\includegraphics[width=0.32\textwidth]{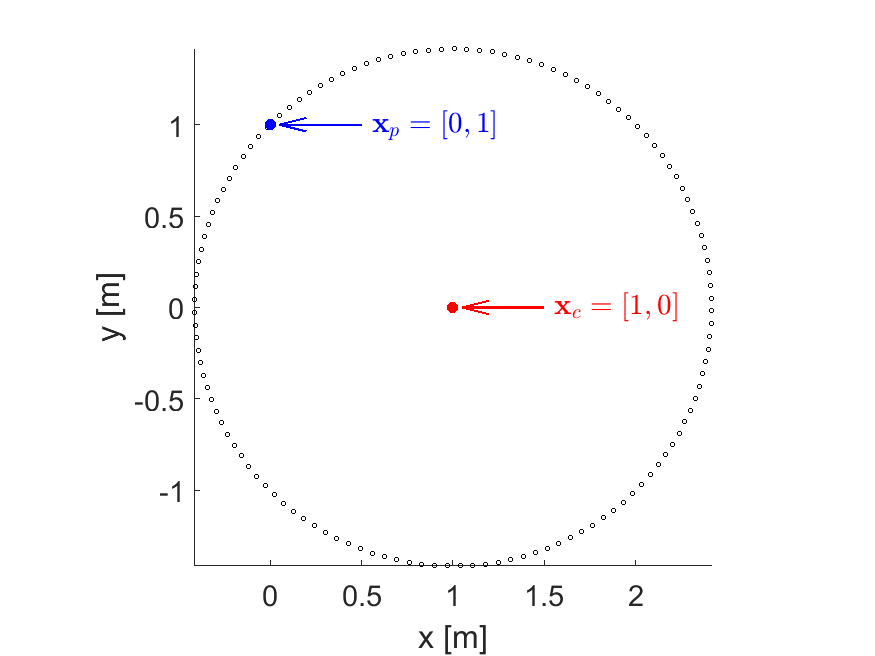} 
\label{fig:2b}  }
\subfigure[]{\includegraphics[width=0.45\textwidth] {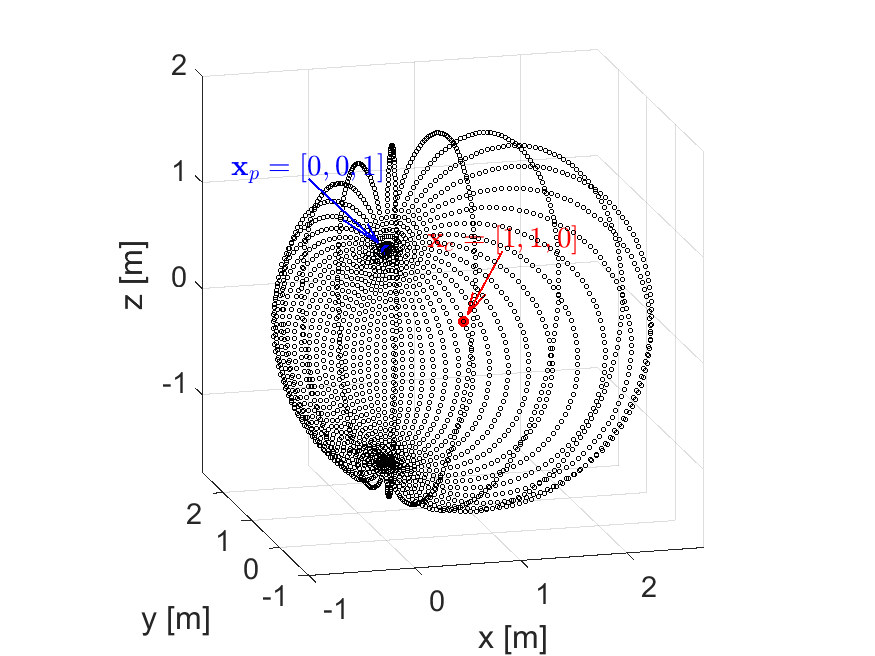} 
\label{fig:2c}  }\\
\subfigure[]{\includegraphics[width=0.45\textwidth]{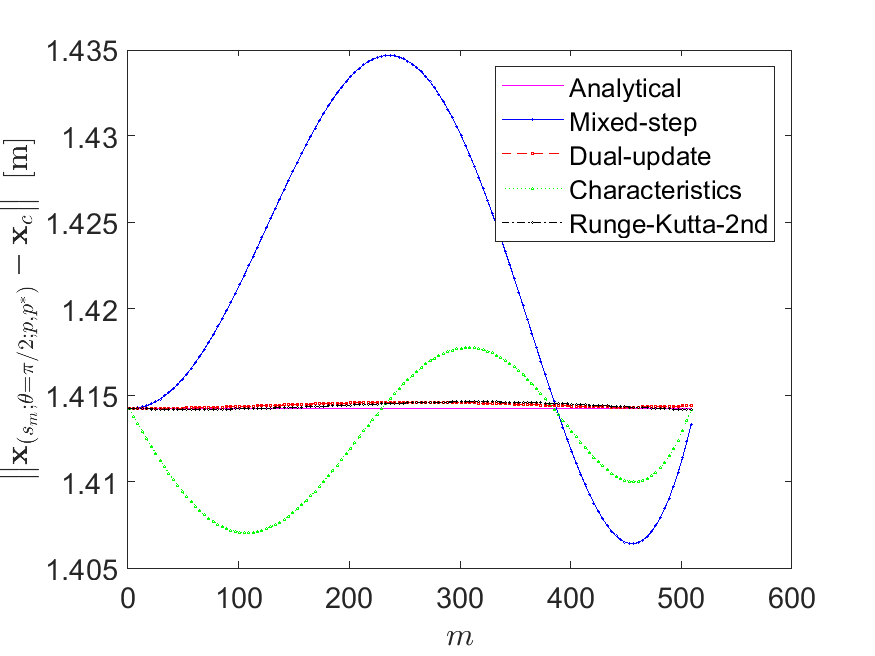}
\label{fig:2d} }
\subfigure[]{\includegraphics[width=0.45\textwidth]{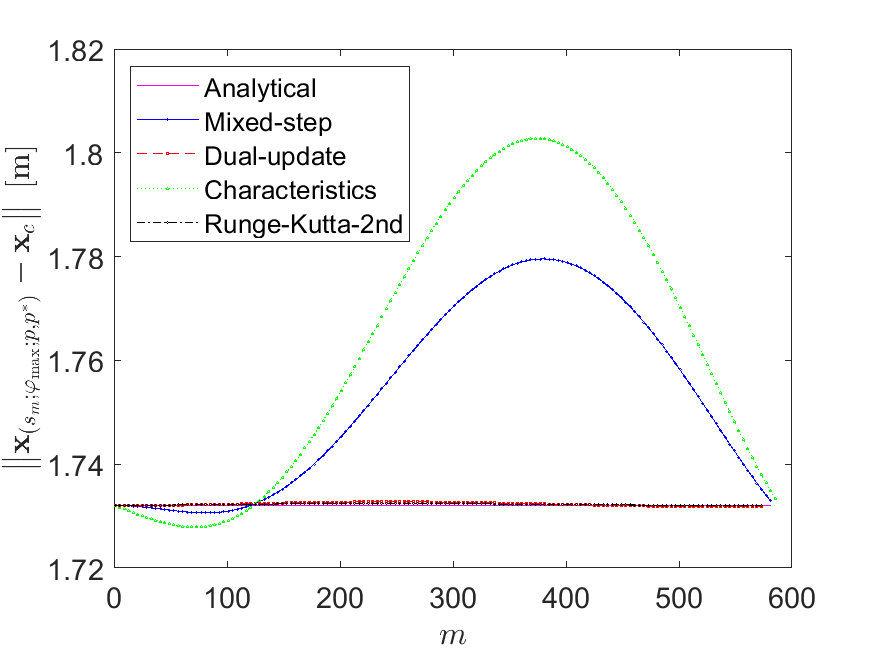}
\label{fig:2e} }
\caption{
(a) Maxwell fish-eye lens. Sampled ray points computed using the second-order Runge--Kutta ray-tracing algorithm: 
(b) 2D case, (c) 3D case. 
Distances of sampled points to the centers of the expected circular paths for four different ray tracing algorithms: 
(d) 2D case, (e) 3D case. 
The ratio $\Delta s / \Delta x = 1$ is used.
}
 \label{fig:2}
\end{figure}

With these settings, the rays are expected to follow circular paths centered at \(\mathbf{x}_c\), returning to their initial positions after a complete $2\pi$ rotation. Figures~\ref{fig:2b} and~\ref{fig:2c} show the ray trajectories for the 2D and 3D experiments, respectively.

Let $\Delta s$ denote the ray spacing. For the trajectories shown in Figures~\ref{fig:2b} and~\ref{fig:2c}, the ray-to-grid spacing was set to unity, i.e., $\Delta s / \Delta x = 1$. For visualization, the points displayed along each ray were sampled at intervals of \(4\Delta s\), i.e., four times the ray step size. Figures~\ref{fig:2d} and~\ref{fig:2e} depict the distances of these sampled points to the expected circular path center for four different ray-tracing algorithms in the 2D and 3D phantoms, respectively. In Figure~\ref{fig:2e} (3D case), the distances correspond to the ray initialized at $\varphi_{\mathrm{max}}$.

For ray point indices \( m \in \{0, \dots, M\} \), the \textit{mean radius deviation} of the circular paths is defined as:

\begin{align}
\mathrm{RE}_{\mathrm{rd}} 
&= \mathrm{mean}_{\varphi} \left[ \frac{1}{M_{(p,p^*)}} 
\sum_{m=1}^{M_{(p,p^*)}} 
\frac{\left| \left\| \mathbf{x}_{(s_m; \varphi; p,p^*)} - \mathbf{x}_c \right\| - R_{\mathrm{True}} \right|}{R_{\mathrm{True}}} \right] \times 100,
\end{align}
where \(p^*\), the last point along the ray, satisfies \(\|\mathbf{x}_p - \mathbf{x}_{p^*}\| < \Delta s\), and \(R_{\mathrm{True}} = \sqrt{\mathrm{dim}}\, a\).
The operator $\mathrm{mean}_{\varphi}$ indicates averaging over the set of rays with different initial azimuthal angles $\varphi$ and is applied only in the 3D case. Throughout this section, we set $a = 1$.

\begin{figure} 
\subfigure[]{\includegraphics[width=0.45\textwidth]{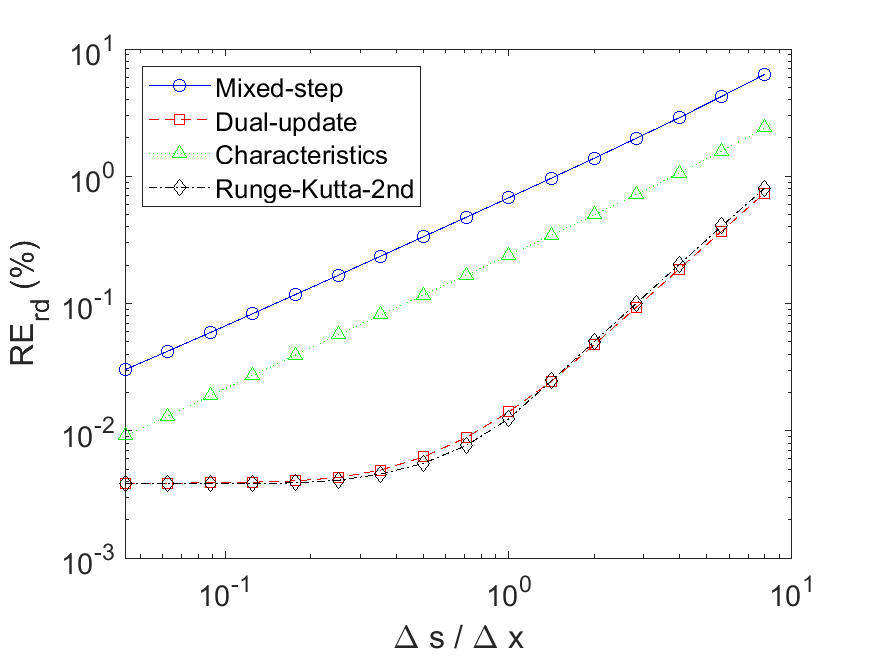}
\label{fig:3a}  }
\subfigure[]{\includegraphics[width=0.45\textwidth]{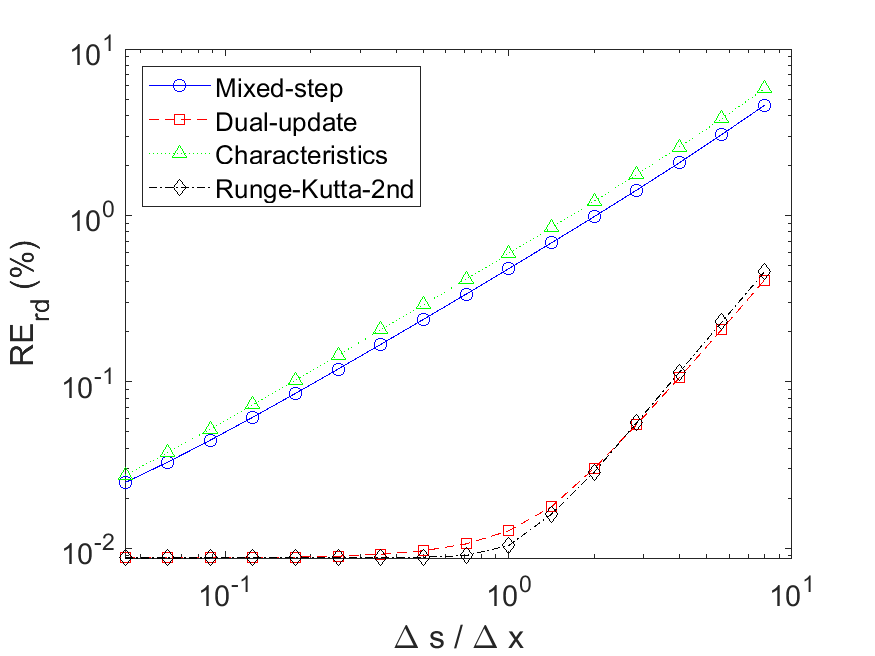}
\label{fig:3b}  }
\caption{
Mean radius deviation (\%) of sampled ray points versus ray-to-grid spacing for four different ray-tracing algorithms. 
}
 \label{fig:3}
\end{figure}

This experiment is repeated for various ray spacings~$\Delta s$. Figures~\ref{fig:3a} and~\ref{fig:3b} show the corresponding $\mathrm{RE}_{\mathrm{rd}}$ values versus the ray-to-grid spacing ratios for the 2D and 3D cases, respectively. As shown in these plots, the Mixed-Step algorithm performs the worst in following the expected analytic trajectories. In contrast, the Dual-Update and Runge–Kutta–2nd algorithms perform the best. However, as the ray-to-grid spacing ratio approaches zero, the interpolation errors begin to dominate, and the contribution of interpolation to the deviation of the ray trajectories becomes larger than the contribution from the ray-tracing algorithms themselves.

\subsection{Mean Acoustic Length Deviation}

The second criterion evaluates how well the traced rays approximate the true acoustic length.  
We define the ray’s initial position as:
\begin{itemize}
    \item[1)] $\mathbf{x}_p = [0, a]$ \quad (\(\mathrm{dim}=2\)),
    \item[2)] $\mathbf{x}_p = [0, 0, a]$ \quad (\(\mathrm{dim}=3\)),
\end{itemize}
and the intended center of the ray’s circular path as:
\begin{itemize}
    \item[1)] $\mathbf{x}_c = [0, 0]$ \quad (\(\mathrm{dim}=2\)),
    \item[2)] $\mathbf{x}_c = [a, a, 0]$ \quad (\(\mathrm{dim}=3\)).
\end{itemize}

\noindent
The initial directions of the rays are set as follows:
\begin{itemize}
    \item[1)] In 2D: \(101\) rays with unit initial direction vectors at angles evenly distributed between \(-\pi/3\) and \(\pi/3\) with respect to the vector $\mathbf{x}_c - \mathbf{x}_p$, i.e., \(\theta \in [-\pi/3, \pi/3] \).
    \item[2)] In 3D: \(21\) rays with directions normal to $\mathbf{x}_c - \mathbf{x}_p$, evenly distributed over a full \(2\pi\) rotation around this axis, i.e., \( \theta = \pi/ 2\) and \(\varphi \in [0, 2\pi) \).
\end{itemize}
\noindent
With these settings, the rays are expected to follow their respective circular paths and be intercepted at:
\begin{itemize}
    \item[1)] In 2D: $\mathbf{x}_{\tilde{p}} = [0, -a]$, after traversing chords tangent to the periphery of circles centered along axes normal to the initial directions.
    \item[2)] In 3D: $\mathbf{x}_{p} = [0, 0, a]$, after completing a full \(2\pi\) rotation.
\end{itemize}

The \emph{mean acoustic length deviation} is defined as:
\begin{align} \label{RE_al}
\mathrm{RE}_{\mathrm{al}} = \mathrm{mean}_{(\theta,\varphi)} \left[ 
\frac{ L_{\left((\theta,\varphi); p, p^+\right)} - L_{\mathrm{True}} }{L_{\mathrm{True}}}
\right] \times 100,
\end{align}
where \(L_{\left((\theta,\varphi); p, p^+\right)}\) is the accumulated acoustic length of the ray connecting points \(p\) to \(p^+\), initialized at varying angles \(\theta\) and \(\varphi\) for the 2D and 3D cases, respectively.  
Also, we set \(p^+ = \tilde{p}\) for the 2D case and \(p^+ = p\) for the 3D case.  
To ensure this connection, the last point is replaced by \(p^+\), yielding a ray spacing \(\Delta s' < \Delta s\) between the two last points.

\noindent
The acoustic length of a ray connecting \(p\) to \(p^+\) is approximated using the trapezoidal rule:
\begin{align}
\begin{split}
L_{((\theta,\varphi); p, p^+)} =& \frac{\Delta s}{2} \, n\!\left(\mathbf{x}_{(s_0; (\theta,\varphi); p, p^+)}\right) 
+ \Delta s \sum_{m=1}^{M_{(p, p^+)}-2} n\!\left(\mathbf{x}_{(s_m;(\theta,\varphi) ; p, p^+)}\right) \\
&+ \frac{\Delta s + \Delta s'}{2} \, n\!\left(\mathbf{x}_{(s_{M-1}; (\theta,\varphi); p, p^+)}\right) 
+ \frac{\Delta s'}{2} \, n\!\left(\mathbf{x}_{(s_M; (\theta,\varphi); p, p^+)}\right).
\end{split}
\end{align}

\noindent
The corresponding analytic acoustic length is:
\begin{align}
L_{\mathrm{True}} = 2a(\mathrm{dim} - 1) \arctan(1) = a(\mathrm{dim} - 1) \frac{\pi}{2},
\end{align}
where \(L_{\mathrm{True}}\) depends on \(\mathrm{dim}\) only due to the distinct interception points chosen for the 2D and 3D experiments.

\begin{figure} 
\centering
\subfigure[]{\includegraphics[width=0.45\textwidth]{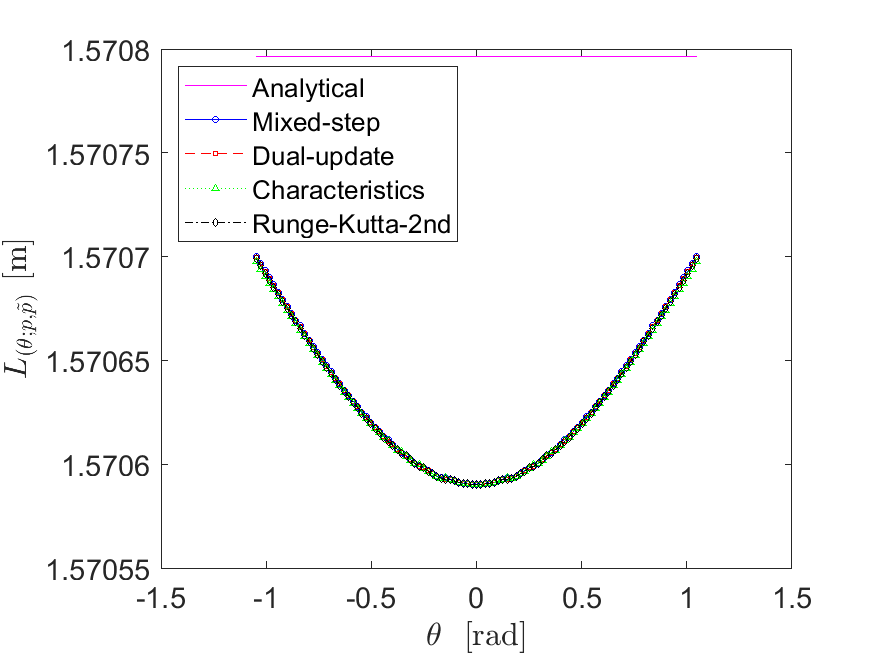}
\label{fig:4a} }
\subfigure[]{\includegraphics[width=0.45\textwidth]{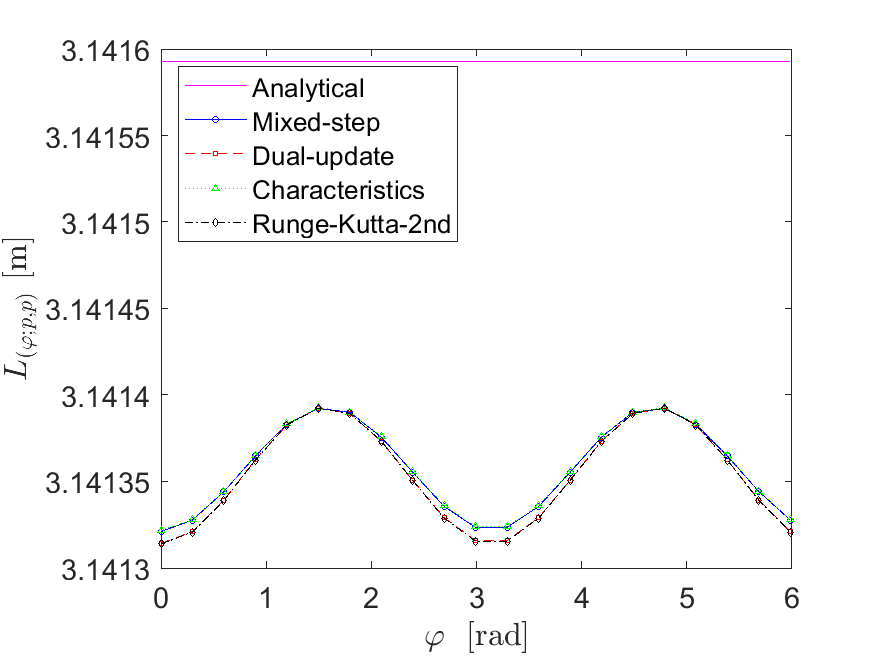}
\label{fig:4b}  }
\caption{
Accumulated acoustic length of rays traced using four ray-tracing algorithms as a function of the initial angles \(\theta\) and \(\varphi\) for 2D and 3D cases, respectively. The ratio \( \Delta s / \Delta x = 1\) is used.
} \label{fig:4}
\end{figure}

\begin{figure} 
\centering
\subfigure[]{\includegraphics[width=0.45\textwidth]{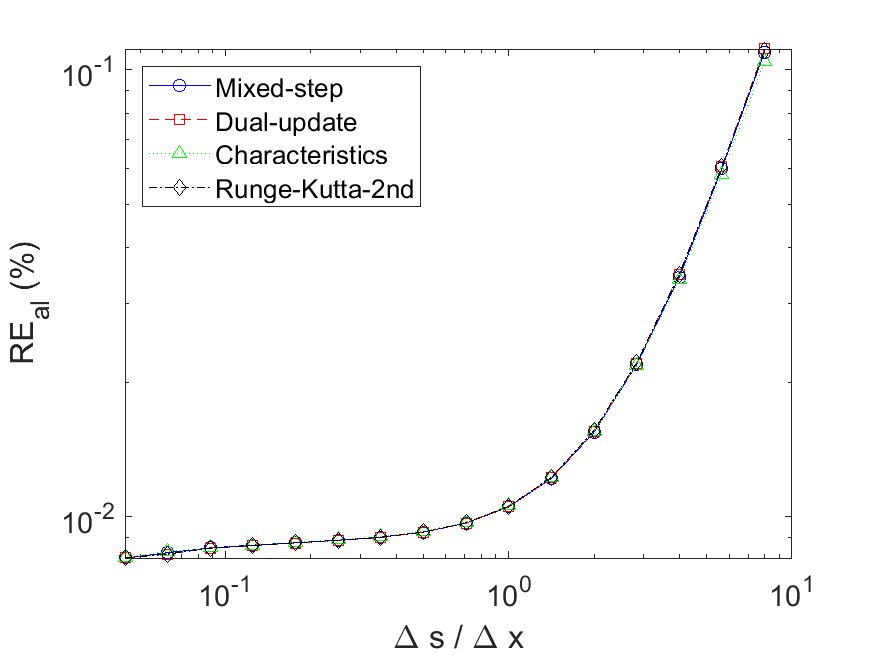}
\label{fig:5a} }
\subfigure[]{\includegraphics[width=0.45\textwidth]{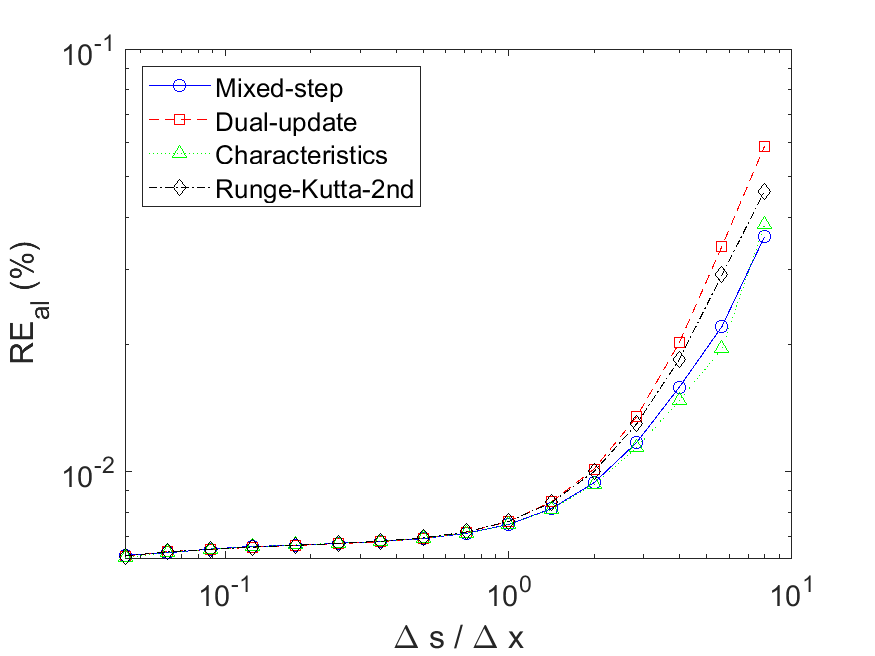}
\label{fig:5b}  }
\caption{
Mean acoustic length deviation (\%) of sampled ray points versus ray-to-grid spacing for four ray-tracing algorithms.
}
 \label{fig:5}
\end{figure}

Figures~\ref{fig:4a} and~\ref{fig:4b} show the acoustic lengths of the rays initialized by varying angles \(\theta \) and \(\varphi \) for the 2D and 3D cases, respectively. Figures~\ref{fig:5a} and~\ref{fig:5b} show the corresponding $\mathrm{RE}_{\mathrm{al}}$ values versus the ray-to-grid spacing for the 2D and 3D cases, respectively.

Since the errors arising from numerical integration along the rays dominated, all four ray-tracing algorithms produced nearly identical accumulated acoustic lengths. This demonstrates that deviations in the ray trajectories contribute negligibly compared to the numerical integration errors in computing the accumulated acoustic lengths.

\section{Discussion and Conclusion}

An open-source \textsc{MATLAB} package was introduced for ultrasound tomography based on two-point ray tracing \cite{Javaherian0}. While the term \emph{ultrasound tomography} encompasses both qualitative and quantitative mapping of acoustic properties from ultrasound data, the primary focus in developing this toolbox was the quantitative reconstruction of the sound speed. Nevertheless, extending it to reconstruct quantitative absorptivity maps or qualitative reflectivity maps is straightforward.  

The toolbox implements two main approaches for reconstructing sound speed images. The first is a time-of-flight (ToF)-based approach, in which a sequence of linearized subproblems is solved, by default, using the SART algorithm~\cite{Anderson}. This yields a low-resolution, low-artefact image of the sound speed. An image reconstructed from the early iterations of the ToF-based approach is then used as the initial guess for the second method: the ray-Born approach \cite{Javaherian-scattering,Javaherian-scattering2}..  

The ray-Born approach combines ray theory, based on the high-frequency approximation, with the Born approximation. From low to high frequencies, the misfit function is linearized, and each linearized subproblem is expressed in terms of forward and backward Green's functions emanating from emitters and receivers, respectively \cite{Javaherian-scattering,Javaherian-scattering2}. Unlike the prototype Born approximation approach, in which Green's functions are computed under the assumption of a homogeneous water medium \cite{Fan,Faucher}, the ray-Born approach uses Green's functions approximated along the ray paths linking emitter–receiver pairs, under the high-frequency approximation \cite{Javaherian-scattering,Javaherian-scattering2}.

The \emph{Numerical Results} section focused on validating the accuracy of the ray-tracing algorithms, which are a core component of our image reconstruction methods. While the second-order Runge--Kutta and Dual-update methods provide more precise path tracing, all four algorithms produce comparable acoustic lengths, especially at the unit ray-to-grid spacing used in image reconstruction. Readers are referred to related works for additional details about this software package~\cite{Javaherian-refraction,Javaherian-scattering,Javaherian-scattering2}.

The image-reconstruction approaches proposed and implemented in this toolbox \cite{Javaherian-refraction, Javaherian-scattering, Javaherian-scattering2} have been applied to in-vitro and in-vivo datasets released by the University of Rochester Medical Center, with the corresponding results presented in \cite{Javaherian-experimental}.

\section*{Data \& Code Availability}
The \textsc{MATLAB} codes supporting the numerical results presented in this study are fully reproducible by running the example script \verb|example_validate_ray_fish_eye_phantom.m| available in the open-source GitHub repository: \\
\url{https://github.com/Ash1362/ray-based-quantitative-ultrasound-tomography}.

\section*{Acknowledgment}
I would like to express my sincere gratitude to Professor Nebojsa Duric and his team at the University of Rochester Medical Center, as well as to Delphinus Medical Technologies, for publicly releasing the in vitro and in vivo datasets that enabled the experimental validation of this toolbox in a parallel related study \cite{Javaherian-experimental}. I also gratefully acknowledge Professor Mohammad Mehrmohammadi, Dr. Rehman Ali, and Mr. Gaofei Jin for their valuable advice and for sharing their experience with experimental ultrasound data.

\section*{Appendix: Equations Governing the Forward Operator}

As described in Section~\ref{sec:software-description}, the forward operator in this study is based on a ray-theoretic model using a high-frequency approximation, solved in the frequency domain \cite{Javaherian-scattering, Javaherian-scattering2}. At each angular frequency \( \omega \), the propagation of acoustic waves produced by an emitter \( e \), modeled as a point source located at \( \bx_e \), is governed by the Helmholtz equation \cite{Javaherian-scattering, Javaherian-scattering2}
\begin{align} \label{eq:wave}
\left[ \tilde{k}(\bx)^2 + \nabla^2 \right] \hat{p}(\omega, \bx, \bx_e) = -s(\omega, \bx_e),
\end{align}
where \( s(\omega, \bx_e) \) is the source term at emitter position \( \bx_e \) and frequency \( \omega \),  
\( \hat{p}(\omega, \bx, \bx_e) \) is the frequency-domain pressure produced by this source, and \( \tilde{\bk} =\bk + \iu \bk_i \) is the complex wavevector, satisfying \cite{Javaherian-scattering, Javaherian-scattering2}
\begin{align}
| \bk | = k, 
\quad \bk_i = \alpha \frac{\bk}{k}  (k_i = \alpha),
\end{align}
where \(\tilde{k} = k + \alpha \iu \) is the wavenumber.
Without loss of generality, we adopt an absorption model based on the frequency power law, given by \( \alpha = \alpha_0 \omega^y \), where \( \alpha_0 \) is the spatially varying absorption coefficient (units: \( \mathrm{Np} \, (\mathrm{rad}/\mathrm{s})^{-y} \, \mathrm{m}^{-1} \)), and \( y \) is a non-integer exponent, typically in the range \( 1 < y \leq 1.5 \) for soft tissues \cite{Szabo}.
\noindent
For the point source \( s(\omega, \bx_e) \), the analytic solution of the Helmholtz equation is
\begin{align} \label{eq:wave-solution}
    \hat{p}(\omega, \bx, \bx_e) = g(\omega, \bx, \bx_e) \, s(\omega, \bx_e),
\end{align}
where \( g(\omega, \bx, \bx') \) is the Green’s function, expressed as
\begin{align}
    g(\omega, \bx, \bx') = A(\omega, \bx, \bx') \, 
    \exp\!\left[ \iu \phi (\omega, \bx, \bx') \right].
\end{align}
\noindent
Here, \(\phi\) is the phase, satisfying
\begin{align}
\bk = \nabla \phi,
\end{align}
and \( A \) denotes the amplitude decay, factorized as
\begin{align}
    A = A_{\mathrm{geom}} \, A_{\mathrm{abs}},
\end{align}
where \( A_{\mathrm{geom}} \) and \( A_{\mathrm{abs}} \) represent the effects of geometric spreading and physical absorption, respectively \cite{Javaherian-scattering, Javaherian-scattering2}.
\noindent
Substituting Eq.~\eqref{eq:wave-solution} into Eq.~\eqref{eq:wave} and applying a high-frequency approximation yields:

1. The \emph{Eikonal equation}:
\begin{align} \label{eq:eikonal}
    \nabla \phi \cdot \nabla \phi = k^2,
\end{align}

2. The \emph{Transport equation}:
\begin{align} \label{eq:transport}
    \nabla \cdot \left[ A_{\mathrm{geom}}^2 \, \nabla \phi \right] = 0.
\end{align}

While \(\phi\) and \( A_{\mathrm{abs}} \) are computed along the rays linking emitter \( e \) to each receiver \( r \) via the Eikonal equation~\eqref{eq:eikonal} \cite{Javaherian-scattering, Javaherian-scattering2},  
\( A_{\mathrm{geom}} \) is obtained from the Jacobian of the ray, evaluated by tracing a paraxial ray in parallel to each reference ray \cite{Javaherian-scattering2}.

The minimization of the linearized forward operator is formulated using forward and backward Green’s functions \cite{Javaherian-scattering} and their variants incorporating reciprocal phase and amplitude \cite{Javaherian-scattering2}. The backward rays, initialized at the receiver locations, are not explicitly traced; instead, they are obtained by reversing the positions and accumulated parameters along the forward (linked) rays initialized at the emitters. Similar to the forward rays, computing the Jacobian of the backward rays—used to approximate the geometrical spreading of the backward Green’s functions—requires tracing a paraxial ray along each reversed linked ray.

In the time-of-flight formulation, dispersion effects are neglected, and the Eikonal equation reduces to:
\begin{align}
   \nabla L \cdot \nabla L = n^2 
   \quad \left( \nabla T \cdot \nabla T = c^{-2} \right),
\end{align}
where \( L(\bx, \bx_e) \) is the acoustic length and \( T(\bx, \bx_e) \) is the travel time between the source point \( \bx_e \) and position \( \bx \) \cite{Javaherian-refraction}.


\begin{thebibliography}{9}

\bibitem{Li-tof-imag}
C. Li, N. Duric, P. Littrup and L. Huang, “In-vivo breast sound speed imaging with ultrasound computed tomography,” \textit{Ultrasound in Med. \& Biol.}, vol. 35, no. 10, pp. 16151628, 2009.

\bibitem{Gemmeke}
H. Gemmeke, T. Hopp, M. Zapf, C. Kaiser, N.V. Ruiter, “3D Ultrasound Computer Tomography: Hardware Setup, Reconstruction Methods and First Clinical Results,” \textit{NUCL INSTRUM METH A}, vol. 873, 2017, pp. 59-65, 2017.

\bibitem{Synnevag}
J. F. Synnevag, A. Austeng and S. Holm, “Adaptive beamforming applied to medical ultrasound imaging,” \textit{IEEE T-UFFC}, Vol. 54, no. 8, pp. 1606-1613, 2007. DOI: 10.1109/TUFFC.2007.431.

\bibitem{Wiskin-nonlinear}
J. W. Wiskin, D. T. Borup, E. Iuanow, J. Klock and M. W. Lenox, “3-D Nonlinear Acoustic Inverse Scattering: Algorithm and Quantitative Results,” \textit{IEEE T-UFFC}, vol. 64, no. 8, pp. 1161-1174, Aug. 2017.

\bibitem{Javaherian0}
A. Javaherian, “r-Wave: An open-source MATLAB package for quantitative ultrasound tomography via ray-
Born inversion with in vitro and in vivo validation,” GitHub, 4 November 2022, available at: \url{https://github.com/Ash1362/ray-based-quantitative-ultrasound-tomography/}. (Accessed: 4 November 2022).

\bibitem{Javaherian00}
A. Javaherian, “Transmission ultrasound data simulated using the k-wave toolbox as a benchmark for biomedical quantitative ultrasound tomography using a ray approximation to Green’s function,” available at: \url{https://zenodo.org/records/8330926}, September 2023. (Accessed: 9 September 2023).

\bibitem{Cerveny}
V. Červený, “Seismic ray theory, Cambridge University Press.” 2001.


\bibitem{Javaherian-refraction}
A.~Javaherian et al., “Refraction-corrected ray-based inversion for three-dimensional ultrasound tomography of the breast,” \textit{ Inverse Problems}, vol. 36, 125010, 2020.

\bibitem{Javaherian-scattering}
A.~Javaherian and B. Cox, “Ray-based inversion accounting for scattering for biomedical ultrasound tomography,” \textit{Inverse Problems}, vol. 37, no.11, 115003, 2021.

\bibitem{Javaherian-scattering2}
A.~Javaherian, “Hessian-free ray-Born inversion for high-resolution quantitative ultrasound tomography,” \url{https://arxiv.org/abs/2211.00316}, 2023.

\bibitem{Javaherian-experimental}
A.~Javaherian, “The First In Vitro and In Vivo Validation of the Hessian-Free Ray-Born Inversion for Quantitative Ultrasound Tomography,”
\url{https://arxiv.org/abs/2511.16633}, 2025.


\bibitem{Li-tof-pick}
C. Li, L. Huang, N. Duric, H. Zhang and C. Rowe, “An improved automatic time-of-flight picker for medical ultrasound tomography,” \textit{Ultrasonics}, vol. 49, pp. 61-72, 2009.

\bibitem{Qu}
X. Qu, C. Ren, G. Yan, D. Zheng, W. Tang, S. Wang, H. Lin, J. Zhang and J. Jiang,
“Deep-Learning-Based Ultrasound Sound-Speed Tomography Reconstruction with Tikhonov Pseudo-Inverse Priori,” \textit{Ultrasound in Medicine \& Biology}, Vol. 48, no. 10, pp. 2079-2094, 2022.

\bibitem{Ceccato}
R. C. Ceccato, A. V. Pigatto, R. C. Aster, C. -N. Pai, J. L. Mueller and S. S. Furuie, “Time of Flight Transmission Mode Ultrasound Computed Tomography With Expected Gradient and Boundary Optimization,” \textit{IEEE Trans Biomed Eng}, vol. 72, no. 9, pp. 2720-2731, Sept. 2025.

\bibitem{Huang}
Y. Huang, Y. Zeng, S. Cui, C. Liu and X. Cai, “Geometric a priori informed bent-ray tracing for accelerated sound speed imaging in ultrasound computed tomography,” \textit{Ultrasonics}, Vol. 151, 107595, 2025.





\bibitem{Dapp}
R. Dapp, M. Zapf and N.V. Ruiter, “Geometry Independent Speed of Sound  Reconstruction for 3D USCT Using A priori Information,” Proc. IEEE UFFC Symp. (2011) 1403–1406.




\bibitem{Ali-2024}
R. Ali et al., “2-D Slicewise Waveform Inversion of Sound Speed and Acoustic Attenuation for Ring Array Ultrasound Tomography Based on a Block LU Solver,” \textit{IEEE Trans. Med. Imaging}, vol. 43, no. 8, pp. 2988-3000, Aug. 2024. 

\bibitem{Robins}
T. C. Robins, C. Cueto, J. Cudeiro, O. Bates, O. C. Agudo, G. Strong, L. Guasch, M. Warner, M. Tang, “Dual-Probe Transcranial Full-Waveform Inversion: A Brain Phantom Feasibility Study,” \textit{Ultrasound in Medicine \& Biology}, Vol. 49, no. 10, pp. 2302-2315, 2023.

\bibitem{Simonetti}
F. Simonaetti, L. Huang, N. Duric and O. Rama, “Imaging beyond the Born approximation: An experimental investigation with an ultrasonic ring array,” \textit{Phys. Rev. E}. vol. 76, pp. 036601 2007.

\bibitem{Huthwaite}
P. Huthwaite and F. Simonetti, “High-resolution imaging without iteration: a fast and robust method for breast ultrasound tomography,” \textit{J. Acoust. Soc. Am.}, vol. 130, no. 3, pp. 1721-34, 2011.





\bibitem{Thierry}
P. Thierry, S. Operto, and G. Lambare, “Fast 2-D ray+Born migration/inversion in complex media,” \textit{Geophysics}, Vol. 64, No. 1, pp. 162–181, 1999.

\bibitem{Anderson}
A. H. Anderson and A.C. Kak, “Simultaneous algebraic reconstruction techniques (SART): A superior implementation of the ART algorithm,” \textit{Ultrasonic Imaging}, vol. 6, no. 1, pp. 81-94, 1984.


\bibitem{Klimes}
L Klimes, “Grid travel-time tracing: second-order method for the first arrivals in smooth media,” \textit{Pure and Applied Geophysics}, vol. 148, nos. 3/4, 1996.


\bibitem{Moser}
T. J. Moser, G. Nolet and R. Snieder, “Ray bending revisited,” \textit{Bulletin of the Seismological Society of America}, vol. 82, no. 1, pp. 259-288, 1992.

\bibitem{Rawlinson}
N. Rawlinson, G. A. Houseman and C. D. N. Collins, “Inversion of seismic refraction and wide-angle reflection travel-time for three-dimensional layered crustal structure,” \textit{Geophys. J. Int.} vol. 145, pp. 381400, 2001.

\bibitem{Anderson2}
A.H. Anderson  and A.C. Kak, “Digital ray tracing in two-dimensional refractive fields,” \textit{J. Acoust. Soc. Am.}, vol. 72, pp. 1593–606, 1982.


\bibitem{Fan}
Y. Fan and L. Ying, “Solving inverse wave scattering with deep learning,” \textit{Annals of Mathematical Sciences and Applications}, Vol. 7, no. 1, pp. 23–48, 2022.

\bibitem{Faucher}
F. Faucher, C. Kirisits, M. Quellmalz, O. Scherzer and E. Setterqvist, Diffraction Tomography, Fourier Reconstruction, and Full Waveform Inversion, In Handbook of Mathematical Models and Algorithms in Computer Vision and Imaging, 2022, doi. 10.1007/978-3-030-03009-4.

\bibitem{Szabo}
T. L. Szabo, “Time domain wave equations for lossy media obeying a frequency power law,” \textit{J. Acoust. Soc. Am.}, vol. 96, pp. 491–500, 1994.

\end{thebibliography}
\end{document}